\documentclass[10pt]{article}

\usepackage{amsfonts}
\usepackage{amssymb}
\usepackage{amsmath}
\usepackage{indentfirst}
\usepackage{hyperref}
\usepackage{breakurl}
\usepackage[T2A]{fontenc}
\usepackage[utf8]{inputenc}
\usepackage[russian,english]{babel}
\usepackage{indentfirst}
\usepackage[usenames]{xcolor}
\usepackage{afterpage}
\usepackage{float}
\usepackage{hhline}
\usepackage{everypage}

\usepackage{ifpdf}
\ifpdf
   \usepackage{pdfpages}
\else
   \RequirePackage[dvips]{graphicx}
\fi

\newcommand{\bigM}{000}
\newcommand{\bigN}{000}

\renewcommand{\Re}{{\mathrm{Re}}}

\newcommand{\myskip}[1]{}

\newcommand{\myi}{{\mathrm i}}

\newcommand{\myd}{{\mathrm d}}
\ifpdf
  \newcommand{\gra}{png}
\else
  \newcommand{\gra}{eps}
\fi

\newcommand{\glebdone}[1]{}

\newcommand{\yuridone}[1]{}

\title {{Approximation
of Riemann's zeta function by finite Dirichlet series: multiprecision numerical approach}}

\author{Gleb Beliakov
\\
School of Information Technology, Deakin University, \\221 Burwood Hwy, Burwood
3125, Australia \\ \url{gleb@deakin.edu.au}\\ and \\
Yuri Matiyasevich
\\
St.Petersburg Department\\ of Steklov Institute of  Mathematics\\
of Russian Academy of Sciences,  Russia\\ \url{yumat@pdmi.ras.ru}
}
\begin{document}
\maketitle

\ \\[3mm]

\abstract{
The finite Dirichlet series from the title are defined by the condition that they
vanish at as many initial zeroes of the zeta function as possible.
It turned out that such series can produce extremely good
approximations to  the values of Riemann's zeta function inside the
critical strip. In addition,  the coefficients of these series
have remarkable number-theoretical
properties discovered in large scale high accuracy numerical experiments.

So far no theoretical explanation to the
observed  phenomena was found.

}

\section{Introduction}

One of the most important open problems in Number Theory is the famous \emph{Riemann Hypothesis}
stated
in \cite{Riemann1859}.
 At the turn of the century, it was included by David Hilbert as part of  his 8th problem, one among 23 most important, in his opinion,  problems \cite{hilbert} left open for the coming 20th century. The Riemann Hypothesis resisted all numerous attempts to (dis)proof it and was recognized by the Clay Institute as one of the 7 \emph{Millennium problems} \cite{millenium}.

The Riemann Hypothesis, RH for short, is a statement about complex zeroes of \emph{Riemann's  zeta function}. This function can be defined via Dirichlet series
\begin{equation}\zeta(s)=\sum_{n=1}^\infty n^{-s}.
\label{zeta}\end{equation}
This series converges only for $\Re(s)>1$ but the function can be analytically continued to the
whole complex plane with the exception of the point $s=1$ which is its only pole.

The zeta function for real $s$ was studied already by  Leonhard  Euler. In particular, he gave in \cite{euler} another definition of the function via a product, namely,
\begin{eqnarray}
\sum_{n=1}^\infty\frac{1}{n^{s}}&=
&\prod_{p\ \mathrm{prime}}\left(1+\frac{1}{{p^{s}}}+\frac{1}{{p^{2s}}}+\dots\right).
\label{eulerproduct}\end{eqnarray}
This equality
can be viewed as an analytic form of the \emph{\em Fundamental Theorem of Arithmetic}
stating that every natural number has a unique factorization into the product of
powers of primes---just expand the right hand side in \eqref{eulerproduct} and get its   left hand side.

The fact that \emph{Euler product},  the right hand side of \eqref{eulerproduct},  is taken over prime numbers, explains the role played by the zeta function in the study of these numbers. In particular, Euler proved anew  the infinitude of prime numbers, and the beauty of his proof can rival that of the original proof given by Euclid: \emph{if the number of primes were finite, then for $s=1$
the divergent harmonic series, that is, the left hand side of
 \eqref{eulerproduct}, would
 have finite value equal to the right hand side of \eqref{eulerproduct}}.

 Bernhard Riemann
 went further, he showed that the zeta function can be used for the study of  the growth of
 the \emph{prime counting function}
$\pi(x)$ equal to
the number of primes not exceeding $x$. This is a step function
having a jump of size $1$ at each prime number.
\myskip{Like the function $\psi(x)$, the
function $\pi(x)$ has a jump
 at each prime, but only of size 1,
  subtler problems of the distribution of primes. skip}

 In a more transparent way  the relationship between
the  zeroes of the zeta function and distribution of prime numbers
can be
be expressed in terms
  of another step function, $\psi(x)$,  defined  by Pafnutij Chebyshev
in \cite{Chebyshev1852}
as
\begin{eqnarray}
 \psi(x)&=&\ln(\mathrm{LCM}(1,2,...,\lfloor x\rfloor)).
\label{psi}\end{eqnarray}
Similar to $\pi(x)$, this  function also
  has a jump at each prime $p$ but now of increasing   size  $\ln( p)$,  and besides it has a jump of the same size at every power of $p$ as well.
 Hans Carl Friedrich von Mangoldt \cite{Mangoldt1895} proved that
\emph{for non-integer $x$ greater than~$1$}
\begin{eqnarray}
 \psi(x)&=&
 x-\sum_{ \zeta(\rho)=0}
 \frac{x^{\rho}}{\rho}
 -\ln(2\pi).
\label{mangoldt}\end{eqnarray}

According to \eqref{mangoldt}, the growth of the difference $\psi(x)-x$
depends on the real parts
of the zeros of the zeta function.
 Already Euler knew that this function vanishes at
  negative even integers, and they are nowadays called the \emph{trivial zeroes.} Riemann proved that they are the only real zeroes of the zeta
  function and that all other, \emph{non-trivial zeroes}  lie inside the so-called \emph{critical strip} $0\le\Re(s)\le 1$.

   Riemann's Hypothesis predicts that in fact
{the non-trivial zeros
lie on the
\emph{critical line}
$\Re(s)=\frac{1}{2}$.}
In terms of Chebyshev's function RH can be restated
 as
\begin{equation}
\psi(x)=x+O(x^{\frac{1}{2}}\ln^2(x))
\label{RHpsi}\end{equation}
and in terms of the function $\pi(x)$  as
\begin{equation}
\pi(x)=\int^x \frac{\myd t}{\ln(t)}+O(x^{\frac{1}{2}}\ln(x)).
\end{equation}

Many researchers verified
 the validity of  RH
 for  initial  zeroes of the zeta function via  finite computations  giving, nevertheless mathematically rigourously, the exact value   $\frac{1}{2}$ for their real parts.
 The last achievement reported in \cite{gourdon}  tells that this is so for impressive
$10^{13}$ initial (pairs of conjugate) zeroes of the zeta function.

Numerical studies of the zeta function are valuable from the perspective of discovering interesting patterns  in its behaviour, providing  preliminary evidence for undiscovered phenomena,
 and formulating hypotheses that are not obvious from the analytic formulas.
 In this article we followed such an approach, by studying numerically various quantities related to approximation the zeta function by finite Dirichlet series.

The simplest form of such series is just the truncation
\begin{equation}\zeta_N(s)=\sum_{n=1}^N n^{-s}.
\label{trunkzeta}\end{equation}
Paul Turán \cite{Turan48} established that for proving the Riemann Hypothesis it would be sufficient to show that
\begin{equation} \label{eq:turan}
  \sup\{\Re(s):\zeta_N(s)=0\}=1+O(N^{-\frac{1}{2}}).
\end{equation}
However, Hugh Lowell Montgomery \cite{Montgomery83} proved that in fact
\begin{equation}
  \sup\{\Re(s):\zeta_N(s)=0\}=1+\Omega_+\!\!\left({\frac{\ln\ln(N)}{\ln(N)}}\right),
\label{montg}\end{equation}
which implies that \eqref{eq:turan} does not hold, and hence one cannot prove RH in that way.

 Partial sums of  Riemann's
zeta-function were also studied by
Michel Balazard and  Oswaldo Velásquez Castañón
 in \cite{BalazardC}, by
Peter Borwein, Greg Fee, Ron Ferguson, and Alexa Van Der Waall in
 \cite{BorweinFFW},
 by
Steven M. Gonek  and Andrew H. Ledoan in \cite{GonekL},
by  Norman Levinson in \cite{Levinson},
by  Robert Spira in \cite{Spira1966,Spira1968,Spira1972},
 and by Sergej Voronin in \cite{Voronin}.

 In this article we report on numerical studies of coefficients of finite Dirichlet series  that are constructed not by truncating the infinite series \eqref{zeta} but on the basis of a
 few initial  non-trivial zeroes of the zeta function.
  Firstly, we found that such finite Dirichlet series  approximate well many of the subsequent
  non-trivial zeroes and a number of initial trivial zeroes. This finding (originally observed for a slightly different approximation
 in \cite{yumatLeicester}) was quite unexpected.

Secondly, numerical experiments with very high accuracy revealed that these coefficients have very rich fine structure  related to prime numbers.

The article is structured as follows. In Section \ref{objects} we introduce our objects of study. Section \ref{first} describes the initial findings.  Section \ref{strategy} is devoted to technical details of performing the calculations. In Sections \ref{larger}--\ref{fine} we discussed numerically observed phenomena. In Section \ref{other} we briefly present some similar experiments and our plans for new calculations. In Section \ref{conclusion} we summarize our discoveries.

\section{Our objects for examination}\label{objects}

We are to approximate the zeta function  by  finite Dirichlet series having the form
\begin{equation}
\Delta_N(s)=\sum_{n=1}^N \delta_{N,n}n^{-s}
\label{Delta}\end{equation}
with some weight coefficients $\delta_{N,n}$.
 These coefficients will be selected in such a way that the finite series \eqref{Delta} and (the function defined by) infinite series \eqref{zeta} would have $N-1$ common zeroes.

  The non-trivial zeroes come in conjugate pairs:
\begin{equation} \dots =\zeta(\overline{\rho_3})=\zeta(\overline{\rho_2})=\zeta(\overline{\rho_1})=0=
\zeta(\rho_1)=\zeta(\rho_2)=\zeta(\rho_3)=\dots\end{equation}
Assuming that
  they are simple and satisfy RH, we write
\begin{equation}\rho_n=\frac{1}{2}+\myi \gamma_n
\label{rho}\end{equation}
with
\begin{equation}
 \qquad 0<\gamma_1<\gamma_2<\gamma_3\dots
 \end{equation}

 We will always take for $N$ an odd number,
 $N=2M+1$, put $\delta_{N,1}=1$ and determine the remaining coefficients in \eqref{Delta} by the condition
 \begin{equation}
   \Delta_N\left(\frac{1}{2}\pm\myi\gamma_k\right)=0, \quad k=1,\dots, M.
 \label{Deltadef}\end{equation}

  This condition gives
 explicit expressions for the coefficients in \eqref{Delta}, namely,
\begin{equation}
  \delta_{N,n}=\frac{\tilde\delta_{N,n}}{\tilde\delta_{N,1}},
\label{delta}\end{equation}
where
\begin{eqnarray}
 \lefteqn{ \tilde\delta_{N,n}=(-1)^{n+1}\times}\nonumber\\
 && \begin{vmatrix}
  1&1&\dots &
   1&1\\
  \vdots&\vdots &\ddots&
  \vdots&\vdots\\
  (n-1)^{-\overline {\rho_{1}}}&(n-1)^{-{\rho_{1}}}&\dots&
  (n-1)^{-\overline{ \rho_{M}}}&(n-1)^{- \rho_{M}}\\
   (n+1)^{-\overline {\rho_{1}}}&(n+1)^{-{\rho_{1}}}&\dots&
  (n+1)^{-\overline{ \rho_{M}}}&(n+1)^{- \rho_{M}}\\
  \vdots &\vdots &\ddots &\vdots&\vdots\\
  N^{-\overline{ \rho_{1}}}&N^{- \rho_{1}}&\dots&
N^{-\overline{ \rho_{M}}}&N^{- \rho_{M}}
      \end{vmatrix}.
\label{deltatilde}\end{eqnarray}

\section{First observations}\label{first}

Our  interest was to examine numerical values of the determinants \eqref{deltatilde}.
Originally, it was guessed   that with the growth of $N$ the coefficients $\delta_{N,n}$
(defined by \eqref{delta}) will  approach the coefficients from
\eqref{zeta},  that is, for a fixed~$n$
 \begin{equation}
  \delta_{N,n}\underset{N\rightarrow\infty}{\longrightarrow} 1.
\label{guess1}\end{equation}
This guess was based on an expected analogy with the Taylor series.
Namely, if
\begin{equation}
  1+\sum_{n=1}^\infty a_nz^n =\prod_{k=1}^\infty\left(1-\frac{z}{z_k}\right)
\end{equation}
and
\begin{equation}
   1+\sum_{n=1}^N a_{N,n}z^n =\prod_{k=1}^N\left(1-\frac{z}{z_k}\right)
\label{finprod}\end{equation}
then for a fixed $n$
\begin{equation}
  a_{N,n}\underset{N\rightarrow\infty}{\longrightarrow} a_n.
\end{equation}

\begin{centering}
\begin{figure}[h]\begin{centering}\hfill
  \includegraphics[width=.7\textwidth]{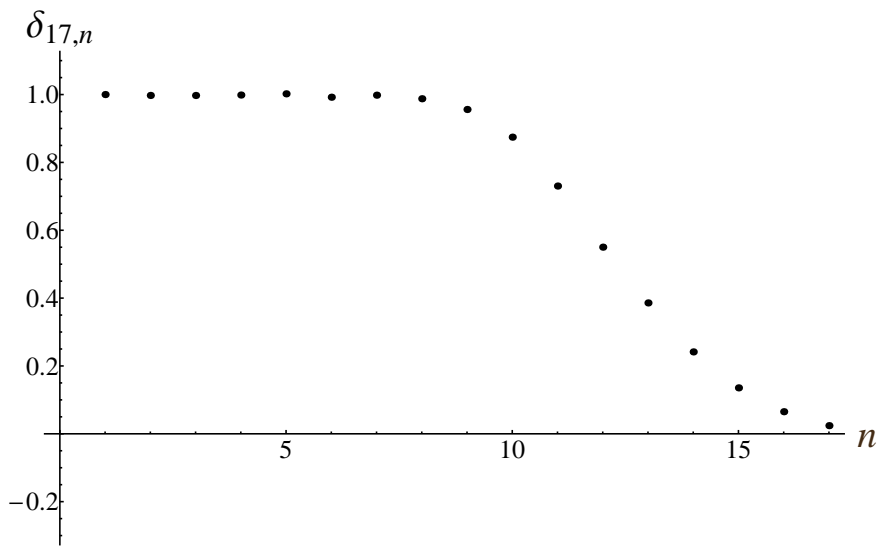}
  \hfill
 \caption{Coefficients $\delta_{17,n}$}
  \label{delta17}\end{centering}
\end{figure}
\end{centering}
Initial calculations seemed to support \eqref{guess1} -- see Figure \ref{delta17}.
This figure justifies our  writing
\begin{equation}\Delta_{17}(s)=\sum_{n=1}^{17}\delta_{17,n}n^{-s}
\leftrightharpoons\sum_{n=1}^{\infty}n^{-s}= \zeta(s)
\label{firstharp}\end{equation}
with the ideograph $\leftrightharpoons$ having here and in the sequel
a very weak sense: \emph{
a few initial coefficients of the two Dirichlet series are
approximately equal}.

\begin{figure}[h]\centering
   \includegraphics[width=.9\textwidth]{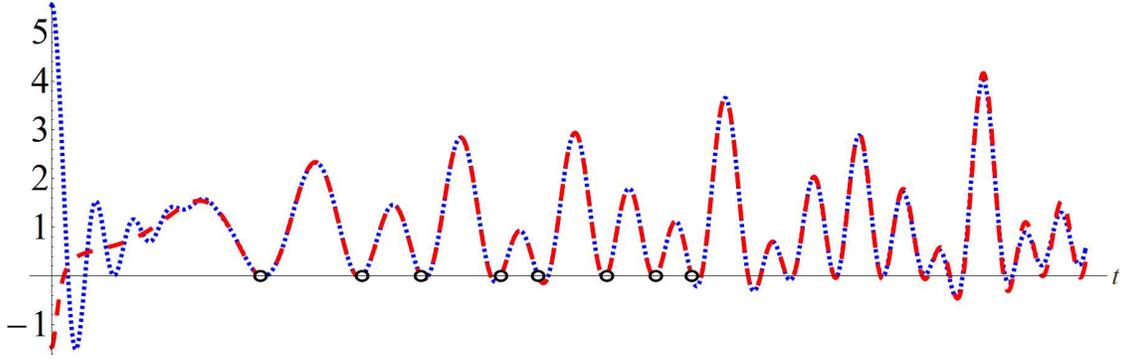}
\caption{Real parts  of \textcolor{red}{$\zeta(\frac{1}{2}+\myi t)$ (dashed)}  and
\textcolor{blue}{$\Delta_{17}(\frac{1}{2}+\myi t)$ (dotted).
}}
\label{re17}\end{figure}

\begin{figure}[h]\centering
  \includegraphics[width=.9\textwidth]{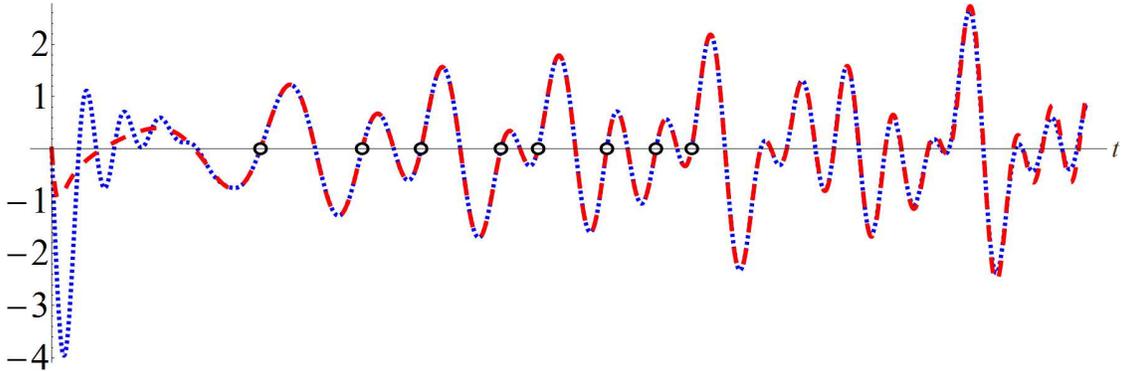}
\caption{Imaginary parts of \textcolor{red}{$\zeta(\frac{1}{2}+\myi t)$ (dashed)} and
\textcolor{blue}{$\Delta_{17}(\frac{1}{2}+\myi t)$ (dotted).
}}\label{im17}
\end{figure}

It turned out that  $\Delta_{17}(s)$ gives a rather good approximation to $\zeta(s)$ on the critical line, see Figures \ref{re17}--\ref{im17}.

\glebdone{Changed captions for Figures in color }


\begin{table}
\begin{center}
\caption{Zeroes of $\Delta_{17}(s)$ nearby zeros of $\zeta(s)$}
\label{appr17}
\vspace{3mm}
\begin{tabular}{rcl}%
    0%
  &%
    =%
  &%
    $\Delta_{17}(\rho_{9}-4.396\dotsc\cdot 10^{-3}+5.711\dotsc\cdot 10^{-3}\myi)$%
  \\ 
    0%
  &%
    =%
  &%
    $\Delta_{17}(\rho_{10}-1.141\dotsc\cdot 10^{-2}-3.345\dotsc\cdot 10^{-3}\myi)$%
  \\ 
    0%
  &%
    =%
  &%
    $\Delta_{17}(\rho_{11}-1.498\dotsc\cdot 10^{-2}+1.762\dotsc\cdot 10^{-3}\myi)$%
  \\ 
    0%
  &%
    =%
  &%
    $\Delta_{17}(\rho_{12}-1.158\dotsc\cdot 10^{-2}+2.264\dotsc\cdot 10^{-2}\myi)$%
  \\ 
    0%
  &%
    =%
  &%
    $\Delta_{17}(\rho_{13}-1.317\dotsc\cdot 10^{-2}+7.545\dotsc\cdot 10^{-2}\myi)$%
  \\ 
    0%
  &%
    =%
  &%
    $\Delta_{17}(\rho_{14}-7.400\dotsc\cdot 10^{-2}-5.559\dotsc\cdot 10^{-4}\myi)$%
  \\ 
    0%
  &%
    =%
  &%
    $\Delta_{17}(\rho_{15}+4.486\dotsc\cdot 10^{-2}+8.379\dotsc\cdot 10^{-2}\myi)$%
  \\ 
\end{tabular}%

 \end{center}
\end{table}

In particular, $\Delta_{17}(s)$ has zeroes close to a few zeta zeroes
following zeroes $\rho_1$, \dots $\rho_8$ used for constructing this finite Dirichlet series, see Table~\ref{appr17}. Nothing similar can happen for Taylor
series -- clearly, the finite product in \eqref{finprod}  contains no information about the
subsequent zeroes $z_{N+1}$, $z_{N+2}$, \dots

The closeness of the values of $\Delta_{17}(s)$ and its zeroes to the values of $\zeta(s)$ and its zeroes is surprising for two reasons:
\begin{itemize}
 \item    the  meaning of the relation
$\leftrightharpoons$ in \eqref{firstharp} is very week;
  \item the infinite series in \eqref{firstharp} diverges on the critical line.
\end{itemize}

\section{Numerical strategies and pitfalls} \label{strategy}

The initial observations prompted  more thorough numerical studies of the determinants \eqref{deltatilde} in order to understand better their behaviour for larger $N$. While experimental numerical values can certainly point to some interesting patterns,  inaccurate experimental results can become false leads, that are due solely to numerical artifacts. For this reason we aimed at providing numerical evidence at a very high precision level, ideally with tight error bounds, as to minimise the likelihood of false leads.

We were aware that calculation of the determinants \eqref{deltatilde} could lead to losses of accuracy, and decided to perform calculations with very high precision  of over ten thousand decimal places. Such an accuracy was achieved by using multiprecision arithmetic, implemented in such packages as GMP \cite{GMPLIB}, Arprec \cite{Arprec} and Arb \cite{arb}. This accuracy allowed us to separate numerical artifacts due to the loss of precision in numerical calculations from some interesting phenomena reported in the subsequent sections.

Let us describe our computational settings.
The values $ \delta_{N,n}$ were computed from $\gamma_1,\ldots,\gamma_M$ by calculating a sequence of determinants ($N=1,2,3,\ldots, 12000$) of a matrix with entries $a_{ij}=i^{-\rho_k}$ for even $j=2k$ and $a_{ij}=i^{-\bar \rho_k}$ for odd $j=2k-1$. The determinants were computed by using a variant of Gauss elimination as reported in \cite{BeliakovM}, in multiprecision arithmetics, using ten thousand decimal places accuracy. The values of $\gamma_k$ were precomputed with twenty thousand decimal places by the authors using Newton-based root finding routine by Fredrik Johansson in his new system Arb \cite{arb}.

These values are available at \cite{MatiyasevichB20000} and more accurate values are at \linebreak \cite{BeliakovM40000}.
 The library GMP \cite{GMPLIB} was used for multiprecision arithmetics, and computations were performed in parallel on an MPI-based cluster involving 168 processes and 400 GB combined RAM, thanks to  VPAC and Monash e-research centre \url{http://www.vpac.org} and  \url{http://www.monash.edu.au/eresearch/}.

The loss of accuracy in Gauss elimination was estimated by computing the same quantities with 20000 decimal places (but for smaller $N$ up to $6000$, due to limitations on computer resources available at the time).
 Thus we had an estimate for the accuracy of the computed coefficients $\delta_{N,n}$. This estimate suggested that the chosen accuracy of 10000 decimal places was in fact warranted due to cancelation errors, and that it was also sufficient for numerical studies up to $N=12001$.

However, this accuracy turns out to be insufficient when we are to examine subtler structure of the coefficients presented in Section \ref{fine} -- cf. Figure \ref{plot10001} with Figure \ref{plot12001},
the lower dots on the latter lie on a horizontal line only because of the lost of accuracy.

Of course, given the computational cost of Gauss elimination and of multiprecision arithmetics, a sufficiently large cluster of processors and combined RAM were needed. The details of our computations are presented in \cite{BeliakovM}. Briefly, we were able to compute the sequence of determinants
$ \delta_{N,n}, n=1,\ldots,N, N=1,2,\ldots, \bar N$ in just one Gauss elimination in $O(\bar N^3)$ time and using $\bar N^2+O(\bar N)$ storage. In our computations we took $\bar N=12001$. All computations were parallelised with particular attention to load balancing,  and computations took seven days on a cluster with 168 processes (Intel  E5-2670 nodes with $48-64$ GB of RAM, connected by 4x QDR Infiniband Interconnect, running CentOS 6 Linux).

\section{Discoveries for larger $N$} \label{larger}

The  multiprecision calculations for large $N$ up to 12001 produced interesting findings.
 Firstly, we found that for large $N$ the series
$\Delta_N$ \eqref{Delta} continues to provide very accurate approximations to the zeros of the zeta
function but doesn't approximate its value at other points any longer; these phenomena  will be discussed  in this Section.
 Secondly, high accuracy experiments allowed us to look ``under a microscope'' at the fine structure of the coefficients $ \delta_{N,n}$. Coefficients which unsuspiciously looked as alternating values $+1$ and $-1$ have
 shown very rich fine structure at the precision level between $10^{-1000}$ and $10^{-10000}$, the structure that unexpectedly revealed \emph{prime numbers}! The patterns revealed are so remarkable and regular that we are convinced they could not be due to numerical artifacts. We detail these observations in  the next Section.
\begin{centering}
\begin{figure}[h]\centering
  \includegraphics[width=.7\textwidth]{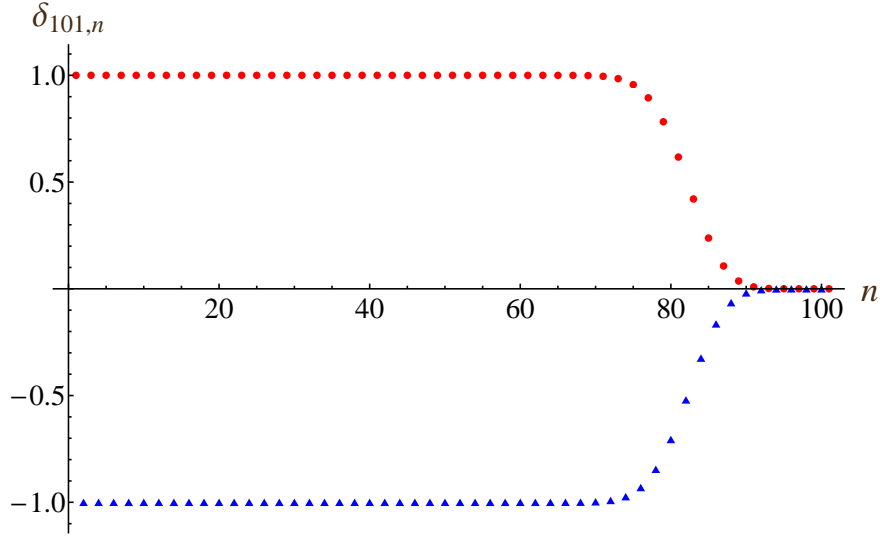}
 \caption{Coefficients $\delta_{101,n}$ \textcolor{red}{(disks} for
  even $n$ and \textcolor{blue}{triangles} for  odd $n$)}
  \label{delta101}


\end{figure}
\end{centering}


\renewcommand{\bigM}{110}
\renewcommand{\bigN}{221}

High accuracy calculations for larger $N$ revealed that most likely the guess \eqref{guess1}
was wrong, and this explains why the values of $\Delta_N(s)$ aren't close to $\zeta(s)$ any longer. Figure~\ref{delta101} exhibits coefficients $\delta_{101,n}$.
Such behaviour is ``typical'' for $N>100$, however, every now and then a kind of ``Gibbs phenomenon''  occurred as illustrated on Figure~\ref{delta\bigN}, or even more
bizarre behaviour as on Figure~\ref{delta233}; presumably, such irregularities would disappear for
$N$ big enough.

\begin{centering}
\begin{figure}[h]\centering
   \includegraphics[width=.7\textwidth]{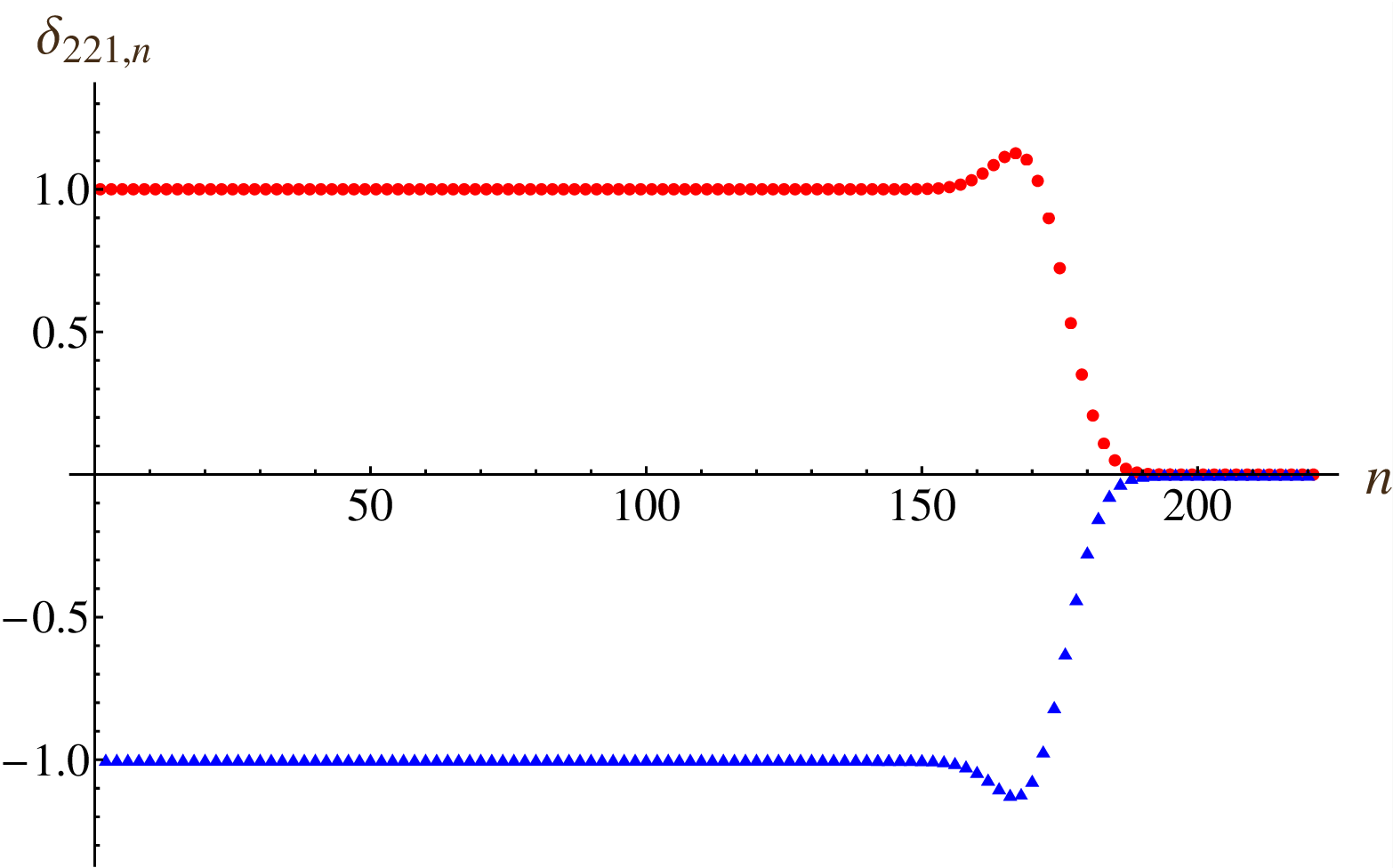}
\caption{Coefficients $\delta_{\bigN,n}$
 (\textcolor{red}{disks} for
  even $n$ and \textcolor{blue}{triangles} for  odd $n$)}
  \label{delta\bigN}
\end{figure}
\end{centering}

\renewcommand{\bigN}{233}

\begin{centering}
\begin{figure}[h]\centering
  \includegraphics[width=.7\textwidth]{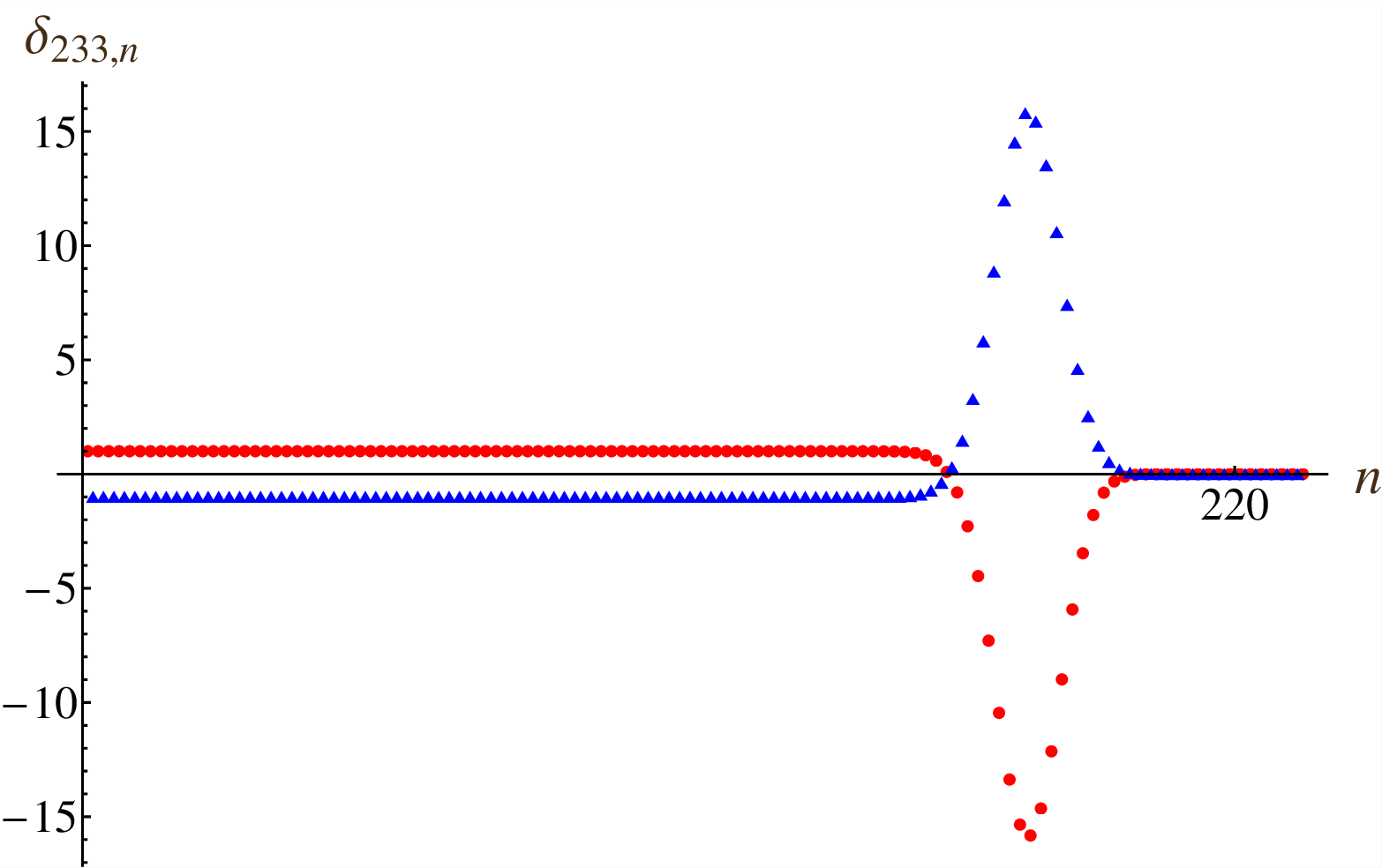}
 \caption{Coefficients $\delta_{\bigN,n}$
 (\textcolor{red}{disks} for
  even $n$ and \textcolor{blue}{triangles} for  odd $n$)}
  \label{delta\bigN}
\end{figure}
\end{centering}

A catalog of $\delta_{N,n}$ for many values of $N$ can be found in
\cite{artless}. Its content suggests that \eqref{guess1} should be replaced by
 \begin{equation}
  \delta_{N,n}\underset{N\rightarrow\infty}{\longrightarrow} (-1)^{n+1}
  \label{altguess1}\end{equation}
and, respectively, for large $N$
\begin{equation}\Delta_{N}(s)=\sum_{n=1}^{N}\delta_{N,n}n^{-s}
\leftrightharpoons\sum_{n=1}^{\infty}(-1)^{n+1}n^{-s}
 = \eta(s)
\label{alt}\end{equation}
where
\begin{equation}\eta(s)=(1-2\cdot 2^{-s}) \zeta(s)
\label{eta}\end{equation}
is the \emph{alternating zeta function}.

\renewcommand{\bigN}{3001}
\begin{table}
\begin{center}
\caption{Zeroes of $\Delta_{\bigN}(s)$ nearby zeroes of $1-2\cdot 2^{-s}$}
\label{Eappr\bigN}
\vspace{3mm}
\begin{tabular}{rcl}%
    0%
  &%
    =%
  &%
    $\Delta_{3001}(1+\frac{2\pi\myi}{\ln(2)}\cdot 50-1.032\dotsc\cdot 10^{-127}+1.020\dotsc\cdot 10^{-127}\myi)$%
  \\ 
    0%
  &%
    =%
  &%
    $\Delta_{3001}(1+\frac{2\pi\myi}{\ln(2)}\cdot 100-2.433\dotsc\cdot 10^{-129}+2.065\dotsc\cdot 10^{-127}\myi)$%
  \\ 
    0%
  &%
    =%
  &%
    $\Delta_{3001}(1+\frac{2\pi\myi}{\ln(2)}\cdot 150+1.032\dotsc\cdot 10^{-127}+1.069\dotsc\cdot 10^{-127}\myi)$%
  \\ 
    0%
  &%
    =%
  &%
    $\Delta_{3001}(1+\frac{2\pi\myi}{\ln(2)}\cdot 200+4.865\dotsc\cdot 10^{-129}+1.146\dotsc\cdot 10^{-130}\myi)$%
  \\ 
    0%
  &%
    =%
  &%
    $\Delta_{3001}(1+\frac{2\pi\myi}{\ln(2)}\cdot 250-1.031\dotsc\cdot 10^{-127}+9.721\dotsc\cdot 10^{-128}\myi)$%
  \\ 
    0%
  &%
    =%
  &%
    $\Delta_{3001}(1+\frac{2\pi\myi}{\ln(2)}\cdot 300-7.294\dotsc\cdot 10^{-129}+2.063\dotsc\cdot 10^{-127}\myi)$%
  \\ 
    0%
  &%
    =%
  &%
    $\Delta_{3001}(1+\frac{2\pi\myi}{\ln(2)}\cdot 350+1.029\dotsc\cdot 10^{-127}+1.117\dotsc\cdot 10^{-127}\myi)$%
  \\ 
    0%
  &%
    =%
  &%
    $\Delta_{3001}(1+\frac{2\pi\myi}{\ln(2)}\cdot 400+9.720\dotsc\cdot 10^{-129}+4.583\dotsc\cdot 10^{-130}\myi)$%
  \\ 
    0%
  &%
    =%
  &%
    $\Delta_{3001}(1+\frac{2\pi\myi}{\ln(2)}\cdot 450-1.027\dotsc\cdot 10^{-127}+9.235\dotsc\cdot 10^{-128}\myi)$%
  \\ 
    0%
  &%
    =%
  &%
    $\Delta_{3001}(1+\frac{2\pi\myi}{\ln(2)}\cdot 500-1.217\dotsc\cdot 10^{-110}+3.892\dotsc\cdot 10^{-111}\myi)$%
  \\ 
    0%
  &%
    =%
  &%
    $\Delta_{3001}(1+\frac{2\pi\myi}{\ln(2)}\cdot 550-1.260\dotsc\cdot 10^{-66}+1.455\dotsc\cdot 10^{-67}\myi)$%
  \\ 
    0%
  &%
    =%
  &%
    $\Delta_{3001}(1+\frac{2\pi\myi}{\ln(2)}\cdot 600-2.580\dotsc\cdot 10^{-36}+2.947\dotsc\cdot 10^{-36}\myi)$%
  \\ 
\end{tabular}%

 \end{center}
\end{table}

\begin{table}
\begin{center}
\caption{Zeroes of $\Delta_{\bigN}(s)$ nearby non-trivial zeroes of $\zeta(s)$}
\label{Zappr\bigN}
\vspace{3mm}
\begin{tabular}{rcl}%
    0%
  &%
    =%
  &%
    $\Delta_{3001}(\rho_{1501}-4.005\dotsc\cdot 10^{-1113}+1.113\dotsc\cdot 10^{-1113}\myi)$%
  \\ 
    0%
  &%
    =%
  &%
    $\Delta_{3001}(\rho_{1601}-5.155\dotsc\cdot 10^{-952}-3.960\dotsc\cdot 10^{-952}\myi)$%
  \\ 
    0%
  &%
    =%
  &%
    $\Delta_{3001}(\rho_{1701}-7.652\dotsc\cdot 10^{-849}+1.788\dotsc\cdot 10^{-848}\myi)$%
  \\ 
    0%
  &%
    =%
  &%
    $\Delta_{3001}(\rho_{1801}+1.966\dotsc\cdot 10^{-766}+3.803\dotsc\cdot 10^{-766}\myi)$%
  \\ 
    0%
  &%
    =%
  &%
    $\Delta_{3001}(\rho_{1901}+1.044\dotsc\cdot 10^{-696}-4.253\dotsc\cdot 10^{-696}\myi)$%
  \\ 
    0%
  &%
    =%
  &%
    $\Delta_{3001}(\rho_{2001}+1.021\dotsc\cdot 10^{-636}-8.184\dotsc\cdot 10^{-636}\myi)$%
  \\ 
    0%
  &%
    =%
  &%
    $\Delta_{3001}(\rho_{2101}-5.402\dotsc\cdot 10^{-582}+8.070\dotsc\cdot 10^{-583}\myi)$%
  \\ 
    0%
  &%
    =%
  &%
    $\Delta_{3001}(\rho_{2201}+9.843\dotsc\cdot 10^{-535}+5.389\dotsc\cdot 10^{-535}\myi)$%
  \\ 
    0%
  &%
    =%
  &%
    $\Delta_{3001}(\rho_{2301}-7.327\dotsc\cdot 10^{-492}-5.590\dotsc\cdot 10^{-491}\myi)$%
  \\ 
    0%
  &%
    =%
  &%
    $\Delta_{3001}(\rho_{2401}+6.471\dotsc\cdot 10^{-452}+8.088\dotsc\cdot 10^{-452}\myi)$%
  \\ 
    0%
  &%
    =%
  &%
    $\Delta_{3001}(\rho_{2501}+1.523\dotsc\cdot 10^{-416}-2.324\dotsc\cdot 10^{-416}\myi)$%
  \\ 
    0%
  &%
    =%
  &%
    $\Delta_{3001}(\rho_{2601}-6.612\dotsc\cdot 10^{-384}-2.011\dotsc\cdot 10^{-384}\myi)$%
  \\ 
    0%
  &%
    =%
  &%
    $\Delta_{3001}(\rho_{2701}+6.698\dotsc\cdot 10^{-354}+3.094\dotsc\cdot 10^{-353}\myi)$%
  \\ 
\end{tabular}%

 \end{center}
\end{table}

\begin{table}
\begin{center}
\caption{Zeroes of $\Delta_{\bigN}(s)$ nearby trivial zeroes of $\zeta(s)$}
\label{Tappr\bigN}
\vspace{3mm}
\begin{tabular}{rcl}%
    0%
  &%
    =%
  &%
    $\Delta_{3001}(-100-8.196\dotsc\cdot 10^{-1220})$%
  \\ 
    0%
  &%
    =%
  &%
    $\Delta_{3001}(-200-4.236\dotsc\cdot 10^{-1017})$%
  \\ 
    0%
  &%
    =%
  &%
    $\Delta_{3001}(-300-4.763\dotsc\cdot 10^{-830})$%
  \\ 
    0%
  &%
    =%
  &%
    $\Delta_{3001}(-400-1.441\dotsc\cdot 10^{-654})$%
  \\ 
    0%
  &%
    =%
  &%
    $\Delta_{3001}(-500-1.187\dotsc\cdot 10^{-488})$%
  \\ 
    0%
  &%
    =%
  &%
    $\Delta_{3001}(-600-4.600\dotsc\cdot 10^{-331})$%
  \\ 
    0%
  &%
    =%
  &%
    $\Delta_{3001}(-700-6.183\dotsc\cdot 10^{-181})$%
  \\ 
    0%
  &%
    =%
  &%
    $\Delta_{3001}(-800-1.648\dotsc\cdot 10^{-51})$%
  \\ 
\end{tabular}%

 \end{center}
\end{table}

Indeed, $\Delta_{\bigN}(s)$ has zeroes
 close to zeroes of the factor $1-2\cdot 2^{-s}$
which have the form $1+\frac{2\pi\myi}{\ln(2)}k,\ k=\pm1,\pm2,\dots $
(see Table \ref{Eappr\bigN}).
It also has zeroes close
 to the non-trivial zeroes $\rho_{1501},\dots$ (see Table \ref{Zappr\bigN}) and, surprisingly, to the trivial zeroes as well (see Table \ref{Tappr\bigN}). In other words, the initial non-trivial zeroes ``feel'' the presence of the pole of the zeta function (canceling
 it  by the factor $(1-2\cdot 2^{-s})$ in \eqref{eta}) and ``know'' about the trivial zeroes not used in
the definitions \eqref{delta}--\eqref{deltatilde}.  Nothing similar can happen for  a meromorphic function
approximated by polynomials with the same zeroes -- they would know nothing about the poles.

\begin{table}
\begin{center}
\caption{Calculation of $\zeta(s)$ via  $\Delta_{\bigN}(s)$   }
\label{valTable\bigN}
\vspace{3mm}
\begin{tabular}{||r|l||}%
  \hhline{|t:==:t|}%
    $\hfill s \hfill$%
  &%
    \rule[-3mm]{0mm}{9mm}{$\hfill \left|\frac{\Delta_{3001}( s \hfill)}{1-2\cdot 2^{-s}}- \zeta(s)\right|\hfill $}%
  \\ 
  \hhline{|:==:|}%
    \rule{0mm}{5mm}$25$%
  &%
    $4.2671\dotsc\cdot 10^{-135}$%
  \\ 
    \rule{0mm}{5mm}$2$%
  &%
    $3.9256\dotsc\cdot 10^{-128}$%
  \\ 
    \rule{0mm}{5mm}$1000\myi$%
  &%
    $4.4184\dotsc\cdot 10^{-128}$%
  \\ 
    \rule{0mm}{5mm}$\frac{1}{2}+10\myi$%
  &%
    $1.0953\dotsc\cdot 10^{-127}$%
  \\ 
    \rule{0mm}{5mm}$-1+100\myi$%
  &%
    $3.6324\dotsc\cdot 10^{-127}$%
  \\ 
    \rule{0mm}{5mm}$-25$%
  &%
    $1.6415\dotsc\cdot 10^{-126}$%
  \\ 
    \rule{0mm}{5mm}$2+1000\myi$%
  &%
    $2.3063\dotsc\cdot 10^{-125}$%
  \\ 
    \rule{0mm}{5mm}$\frac{1}{2}+1000\myi$%
  &%
    $3.9630\dotsc\cdot 10^{-124}$%
  \\ 
    \rule{0mm}{5mm}$-1+1000\myi$%
  &%
    $1.4867\dotsc\cdot 10^{-118}$%
  \\ 
    \rule{0mm}{5mm}$-10+1000\myi$%
  &%
    $8.2377\dotsc\cdot 10^{-103}$%
  \\ 
    \rule{0mm}{5mm}$\frac{1}{2}+5000\myi$%
  &%
    $6.5116\dotsc\cdot 10^{-64}$%
  \\ 
    \rule{0mm}{5mm}$-1+5000\myi$%
  &%
    $2.6548\dotsc\cdot 10^{-59}$%
  \\ 
    \rule{0mm}{5mm}$-10+5000\myi$%
  &%
    $2.5001\dotsc\cdot 10^{-32}$%
  \\ 
  \hhline{|b:==:b|}%
\end{tabular}%

 \end{center}
\end{table}

The values of $\Delta_{\bigN}(s)$ are close to the values
of $\eta(s)$
 for $s$ inside the critical strip and even much to the left of it (see Table \ref{valTable\bigN}). In other words, we have a surprisingly  good approximation to
$\zeta(s)$ of the form
\begin{equation}
  \zeta(s)\approx \frac{\Delta_N(s)}{1-2\cdot 2^{-s}}=
  \frac{\sum_{n=1}^N\delta_{N,n}n^{-s}}{1-2\cdot 2^{-s}}.
\end{equation}
In fact, if we allow more terms in the denominator,
  we can obtain (see \cite{yumatBonn}) much better approximations
\begin{equation}
\zeta(s)\approx
\frac{\sum_{n=1}^N\delta_{N,n}n^{-s}}{\sum_{n=1}^L \mu_{N,n}n^{-s}}
\label{zetafrac}\end{equation}
\glebdone{replaced dividing by division}
for a small value of $L$ where numbers $\mu_{N,n}$ are defined via formal division of
the two Dirichlet series:
\begin{equation}
\frac{\Delta_N(s)}{\zeta(s)} =
  \frac{\sum_{n=1}^N \delta_{N,n}n^{-s}}{\sum_{n=1}^\infty n^{-s}}
 =\sum_{n=1}^\infty \mu_{N,n}n^{-s}.
\label{mudiv}\end{equation}

\section{Fine structure of the coefficients $\delta_{N ,n}$ }\label{fine}

\subsection{Sieve of  Eratosthenes}

\renewcommand{\bigN}{10001}

Clearly, the extreme closeness
of the zeroes and values of the alternating zeta function $\eta(s)$
and that of finite Dirichlet series
$\Delta_N(s)$ is due to the very peculiar values of the coefficients
$\delta_{N,n}$, and  now we are to look at their finer structure
 ``under a microscope''.
To this end  we change to the logarithmic scale -- see Figure~\ref{plot\bigN}.

\begin{centering}
\begin{figure}[H]\centering
   \includegraphics[width=.9\textwidth]{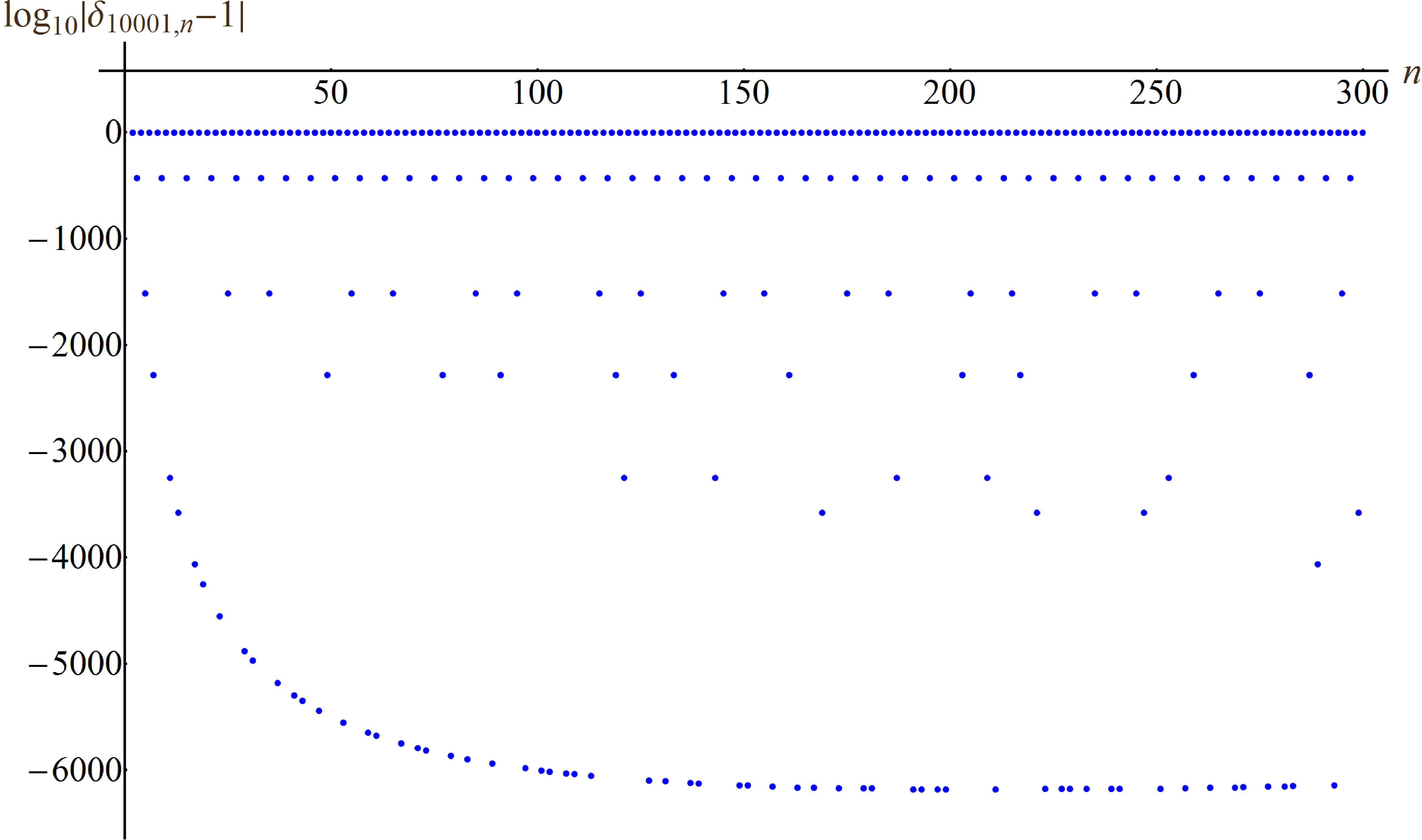}
 \caption{Graph of $\log_{10}|\delta_{\bigN,n}-1|$.}
  \label{plot\bigN}
\end{figure}
\end{centering}

Here we can observe several horizontal rows of dots.
The top row corresponds to even values of $n$ for which $\delta_{\bigN ,n}$
is, according to \eqref{altguess1}, close to~$-1$.
The second row corresponds to odd values of $n$ divisible by~$3$.
The third row corresponds to those values of $n$ that are divisible by $5$
but are relatively prime to $2\cdot 3$.
The fourth row corresponds to those values of $n$ that are divisible by $7$
but are relatively prime to $2\cdot 3\cdot 5$, and so on.
\myskip{The fifth row corresponds to those values of $n$ that are divisible by $11$
but are relatively prime to $2\cdot 3\cdot 5\cdot 7$.
The sixth row corresponds to those values of $n$ that are divisible by $13$
but are relatively prime to $2\cdot 3\cdot 5\cdot 7\cdot 11$.}
The  seventh row, the last one that we can see, contains only two
dots corresponding to $n=17$ and $n=289$.

The remaining dots correspond to prime values of $n$.
So we can say that the initial part of the plot of $\log_{10}|\delta_{\bigN,n}-1|$
represents the \emph{Sieve of  Eratosthenes}.
Respectively, the  horizontal rows
corresponding to the values of $n$
divisible by $2$, by $3$,   \dots but not by the previous primes
will be called
\emph{Eratosthenes levels}.

\begin{centering}
\begin{figure}[H]\centering
\includegraphics[width=.9\textwidth]{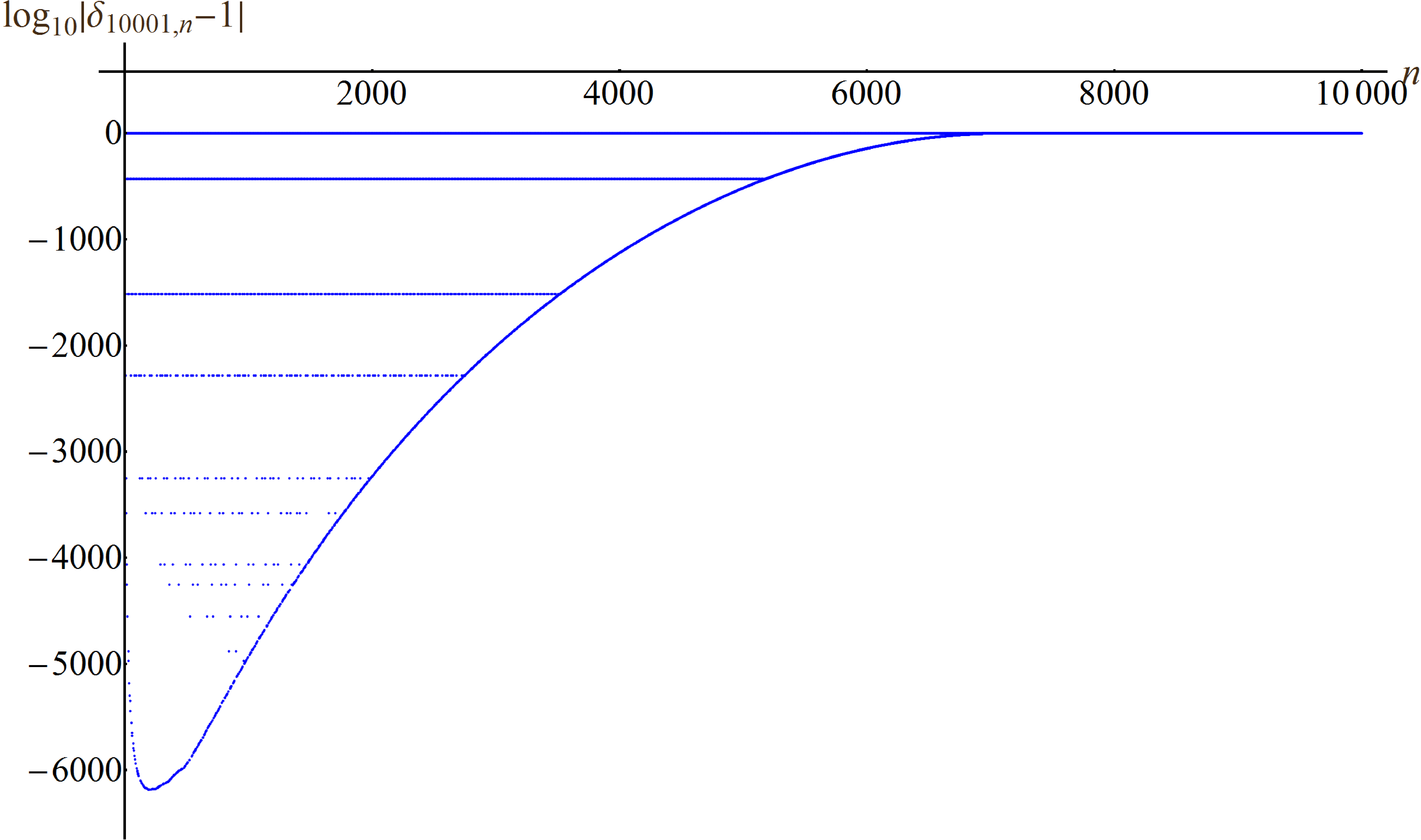}
 \caption{Graph of $\log_{10}|\delta_{\bigN,n}-1|$.}
\label{plot\bigN\bigN}\end{figure}
\end{centering}

Figure~\ref{plot\bigN\bigN}
extends Figure \ref{plot\bigN} up to $n=\bigN$.
 We see that the Eratosthenes levels
break off
when they touch a mysterious ``smooth curve'' of increasing values of
$\log_{10}|\delta_{\bigN,n}-1|$. The larger $N$, the more to the right is the smooth curve.

Figure \ref{plot12001} presents results of our computations for $N=12001$. It again shows the Eratosthenes levels but also gives an impression
of a new phenomenon -- dots corresponding to all primes greater than 80 look like lying on a horizontal
line with the ordinate $-7157$. Actually, this is due to the fact that the calculated values of the
coefficients have only about 7157 correct decimal digits.

\renewcommand{\bigN}{12001}

\begin{centering}
\begin{figure}[H]\centering
\includegraphics[width=.9\textwidth]{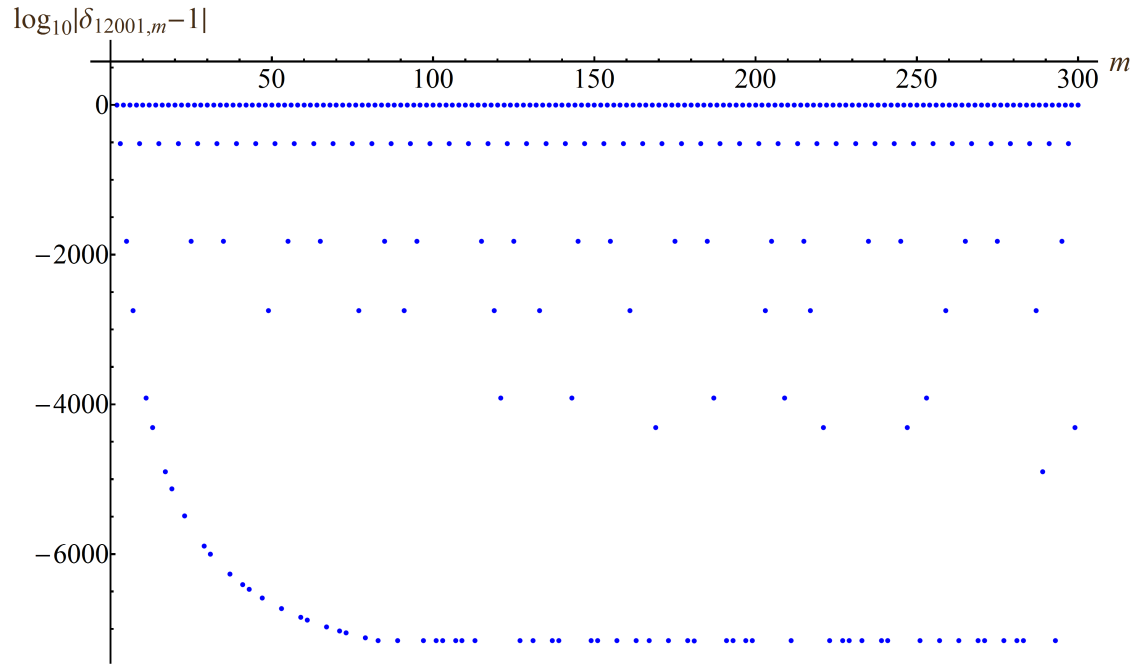}
 \caption{An artifact caused by insufficient accuracy.}
\label{plot\bigN}\end{figure}
\end{centering}

\renewcommand{\bigN}{10001}

\subsection{Fractal structure}

\begin{centering}
\begin{figure}[H]\centering
\includegraphics[width=.9\textwidth]{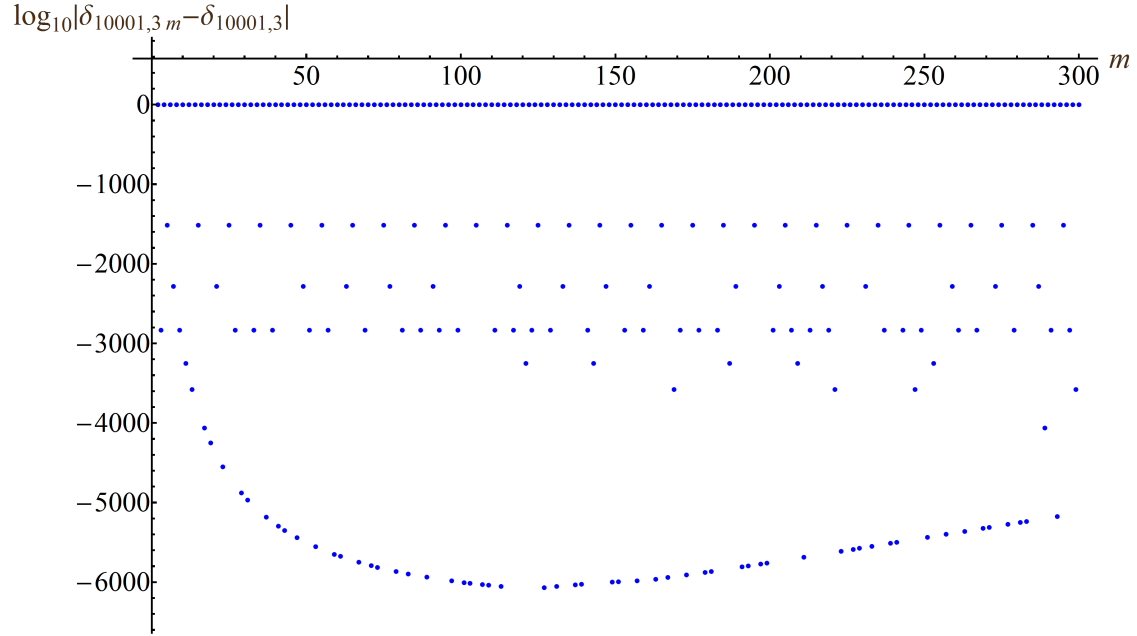}
\caption{Graph of $\log_{10}|\delta_{\bigN,3m}-\delta_{\bigN,3}|$ showing Eratosthenes sublevels.}
  \label{sieve3\bigN}
\end{figure}
\end{centering}

The Eratosthenes levels on Figures \ref{plot\bigN}--\ref{plot\bigN\bigN}
look like lying on straight lines. However, closer examination reveals that each of the levels
in its turn contains sublevels corresponding
to a slightly modified  Sieve of  Eratosthenes.
Figure \ref{sieve3\bigN} shows such sublevels for the main Eratosthenes
level corresponding to  prime $p=3$ in the
case $N=\bigN$. These sublevels correspond to deleting
composite numbers according to their divisibility at first by $2$,
then by $5$, $7$, $3$, $11$, $13$, \dots.

The general rule seems to be as follows. The dots representing $\delta_{N,n}$ for
$n$ from an arithmetical progression $d,2d,\dots,md,\dots$
with $d=2^{k_2}3^{k_3}\dots$
split into Eratosthenes sublevels according to the divisibility of $m$ by
$p_1$,~$p_2$,~\dots where these prime numbers are ordered in such a way that
\begin{equation}
  p_1^{k_{p_1}+1}< p_2^{k_{p_2}+1}<\dots< p_j^{k_{p_j}+1}<\dots\ .
\end{equation}

\section{Related and Further Research}\label{other}

Originally, the second author (\cite{yumatLeicester,yumatBonn} examined the determinants slightly different
from those in \eqref{deltatilde} for which $\Delta^\Gamma_N(s)$, the counterpart of \eqref{Delta}, vanishes at
$2N-2$ zeroes of $\zeta(s)$. This becomes possible thanks to the so called
\emph{functional equation}  established by
Riemann \cite{Riemann1859}. Properties of $\Delta^\Gamma_N(s)$ are similar but not the same as
those of   $\Delta_N(s)$. In particular, the Eratosthenes sieve manifests itself not so
spectacular. On the other hand, $\Delta^\Gamma_N(s)$ allows one to calculate
approximations not only to the zeroes and the values of the zeta function but to its first
derivative as well.

We plan to examine the counterparts of $\delta_{N,n}$ and $\Delta_{N,n}$ for the cases
when zeroes of the zeta function are replaced by zeroes of Dirichlet \mbox{$L$-func}\-tions, as well as to perform computations in interval multiprecision arithmetics using Arb \cite{arb} to obtain rigorous bounds on the resulting values.

This ongoing research can be followed on \cite{artless}.

\section{Conclusion}\label{conclusion}

We performed large scale high accuracy computations of the coefficients of the finite Dirichlet series approximating nontrivial zeros of  Riemann's zeta function. Our aim was to reveal experimentally new relations between these coefficients and various related quantities, such as the zeros of the alternating zeta function. The results of our computations are somewhat unexpected. Firstly they revealed that the finite Dirichlet series also approximates (with high accuracy) other zeros of zeta function (trivial and subsequent non-trivial zeros), not used in computations. Secondly, the coefficients inconspicuously looking as $+1$ and $-1$ have in fact a rich structure related to prime numbers.

We want to underline the necessity for performing computations with very high accuracy, which was crucial in discovering the patterns presented here, that would not be detected otherwise. The calculations performed were costly, of order of 200,000 CPU hours, which were made possible by collaborative work of mathematicians, computer scientists, programmers and support engineers.

Of course, in Number Theory there are many examples of conjectures that were at first substantiated
by calculation for many initial values of the parameters,  but then  were disproved either
theoretically or by finding a numerical counterexample. Nevertheless, we find it highly desirable to  extend our calculations to higher sizes of determinants in order to study subtler
properties of the intriguing numbers  $\delta_{N,n}$. This requires significant computational resources and multi-party collaboration. Our recent experiences with computational aspects of multiprecision calculations are presented in \cite{BeliakovM}.

\section*{Acknowledgements}

The multiprecision values of the zeroes of the zeta function
were first computed  using
 \href{http://www.wolfram.com/mathematica/}{Mathematica}
and
\href{http://www.sagemath.org}{Sage}. The author is grateful to
Oleksandr Pavlyk, special functions developer at Wolfram
Research, for performing part of the calculations, and to
 Dmitrii Pasechnik (NTU Singapore) for hands-on help with
setting up and monitoring the usage of Sage.

More recently  Fredrik Johansson developed a more efficient algorithm
in his new system Arb \cite{arb},
and zeta zeroes were recalculated with 40,000 digits; they are made publicly available
\cite{MatiyasevichB20000,BeliakovM40000} thanks to \href{http://www.monash.edu.au/eresearch/}{Monash e-Research centre} and Multi-modal Australian Sciences Imaging and Visualisation Environment \href{http://www.massive.org.au}{(MASSIVE)}.

  Calculations of zeta zeroes were also performed on computers from
 \href{http://www.grid.am}{ArmNGI}
 (Armenian National Grid Initiative Foundation),
 \href{http://www.newton.ac.uk/}{Isaac Newton Institute for %
Mathematical Sciences}, UK,
 \href{http://lacl.univ-paris12.fr/}{LACL}
  (Laboratoire d'Algorithmique, Complexit\'e et Logique de
  Universit\'e Paris-Est Cr\'eteil),
 \href{http://www.liafa.jussieu.fr/}{LIAFA} (Laboratoire d'Informatique Algorithmique: Fondements et Applications,
 supported jointly by the French National Center for Scientific Research (CNRS)
 and by the University Paris Diderot--Paris 7),
\href{http://www.spiiras.nw.ru/index.php?newlang=english}{SPIIRAS}
(St.Petersburg Institute for Informatics and Automation of RAS),
 \href{http://www.wolfram.com/}{Wolfram
Research}.

The most time-consuming part,  calculation of the sequence of determinants $\delta_{N,n}$, was performed on the ``Chebyshev'' supercomputer
at  Moscow State University
\href{http://hpc.msu.ru/}{Supercomputing Center} and \href{http://www.massive.org.au}{MASSIVE cluster}.

The second author was supported in the framework of the Program of Fundamental Research of the Division of Mathematical Sciences of the Russian Academy of Sciences ``Modern problems of theoretical mathematics''.

\end{document}